\documentclass[11pt]{amsart}
\usepackage{amssymb,latexsym,amsmath,amsthm,enumitem,mathrsfs,geometry}
\usepackage{fancyhdr}
\usepackage{color}
\usepackage{stmaryrd}
\usepackage{hyperref}
\hypersetup{
    colorlinks=true,
    linkcolor=red,
    filecolor=magenta,
    urlcolor=cyan,
    citecolor=green
    }
\usepackage{lineno}
\usepackage{diagbox}
\usepackage{array}
\usepackage{tikz}
\usetikzlibrary{arrows}
\usetikzlibrary{positioning}
\usepackage{float}
\usepackage{bbm}
\usepackage{cleveref}
\usepackage{dsfont}
\usepackage{graphicx}
\usepackage{soul}
\usepackage{verbatim}
\usepackage{comment}
\usepackage{caption} 
\usepackage[symbol*]{footmisc}

\newtheorem{theorem}{Theorem}
\newtheorem{prop}{Proposition}

\theoremstyle{remark}
\newtheorem{remark}{Remark}

\theoremstyle{plain}
\newtheorem*{PCC}{Pair Correlation Conjecture (PCC)}
\newtheorem*{HMH}{Horizontal Multiplicity Hypothesis (HMH)}
\newtheorem*{Rep}{Zero Repulsion}

\definecolor{pink}{rgb}{1,.2,.6}
\definecolor{orange}{rgb}{0.7,0.3,0}
\definecolor{blue}{rgb}{.2,.6,.75}
\definecolor{green}{rgb}{.4,.7,.4}
\definecolor{purple}{RGB}{127,0,255}

\numberwithin{equation}{section}
\allowdisplaybreaks

\title[Pair Correlation Conjecture for simple and critical zeros]{Pair Correlation Conjecture for the zeros of the Riemann zeta-function I: simple and critical zeros}

\author[Goldston]{Daniel A. Goldston}
\address{Department of Mathematics and Statistics, San Jose State University}
\email{daniel.goldston@sjsu.edu}

\author[Lee]{Junghun Lee}
\address{Department of Mathematics Education, Chonnam National University}
\email{junghun@jnu.ac.kr}

\author[Schettler]{Jordan Schettler}
\address{Department of Mathematics and Statistics, San Jos\'{e} State University}
\email{jordan.schettler@sjsu.edu}

\author[Suriajaya]{Ade Irma Suriajaya}
\address{Faculty of Mathematics, Kyushu University}
\email{adeirmasuriajaya@math.kyushu-u.ac.jp}

\dedicatory{``Truth is ever to be found in simplicity, and not in the multiplicity and confusion of things." \\
- Sir Isaac Newton}

\begin{document}

\begin{abstract}
Montgomery in 1973 introduced the Pair Correlation Conjecture (PCC) for zeros of the Riemann zeta-function. He also conjectured that asymptotically 100\% of the zeros are simple. His reasoning to support these two conjectures used the Riemann Hypothesis (RH). Building on Montgomery's approach, Gallagher and Mueller proved in 1978 that PCC under RH implies that 100\% of the zeros are simple. Actually, the method of Gallagher and Mueller does not depend on RH, and thus Montgomery's second simplicity conjecture follows unconditionally from his PCC conjecture. We clarify this result by explicitly not assuming RH and considering PCC as a conjecture only concerning the vertical distribution of zeros. We then show that, for the first time, PCC can also be used to obtain information on the horizontal distribution of zeros. Using Gallagher and Mueller's method and a new idea concerning \lq\lq horizontal multiplicity", we use PCC to prove that asymptotically 100\% of the zeros are not only simple but also on the critical line.
\end{abstract}
\date{\today}

\maketitle


\section{Introduction}

In this paper we do not assume the Riemann Hypothesis (RH) is true or make use of results that depend on RH. The study of pair correlation of zeros of the Riemann zeta-function is concerned with the \emph{vertical distance} between pairs of zeros. If RH is true all the zeros are on the vertical \emph{critical} line with real part $1/2$, and thus the vertical distance between zeros is just the distance between zeros. If RH is false, then the zeros will also have a \emph{horizontal} distribution. The main result of this paper is the following theorem that shows that Montgomery's pair correlation conjecture for the vertical distribution of pairs of zeros dramatically limits the horizontal distribution of zeros.
\begin{theorem} \label{thm1} Assuming the pair correlation conjecture \hyperref[PCC]{\bf PCC}, then asymptotically 100\% of the zeros of $\zeta(s)$ are simple and on the critical line.
\end{theorem}

The proof is based on a method of Gallagher and Mueller \cite{GaMu78} using a result of Fujii \cite{Fu74} which follows work of Selberg \cite{Sel46} who discovered the method for obtaining these results unconditionally without RH. Everything else we need is elementary, and mainly involves removing RH from earlier work where it was not needed.

\section{Zeros of the Riemann zeta-function}

The Riemann zeta-function $\zeta(s)$ has \lq \lq trivial" zeros at $s=-2n$, $n\in \mathbb{N}$, and, letting $s=\sigma+it$ where $\sigma,\, t\in\mathbb{R}$, the remaining zeros are non-real and lie in the \lq \lq critical strip" $0<\sigma < 1$. We denote these \lq \lq non-trivial" zeros by $\rho=\beta+i\gamma$ where $0<\beta<1$ and $\gamma\neq 0$. The Riemann Hypothesis states that every non-trivial zero lies on the \lq \lq critical line" $\sigma = 1/2$, and thus RH is the statement that $\beta = 1/2$ for every non-trivial zero.
Since $\zeta(s)$ is real on the real axis with a simple pole at $s=1$, by the reflection principle any zero $\rho = \beta +i\gamma$ has a reflected zero $\overline{\rho}=\beta-i\gamma$. If $\beta \neq 1/2$ we obtain by the functional equation and the reflection principle four zeros $\rho = \beta +i\gamma$, $\overline{\rho}=\beta-i\gamma$, $1-\rho = 1-\beta-i\gamma$, and $1-\overline{\rho}= 1-\beta+i\gamma$ symmetric with respect to both the real axis and the critical line. 
As usual, we only need to consider zeros with $\gamma >0$ above the real axis which we divide into two disjoint sets: 1) zeros on the critical line with $\beta = 1/2$, and 2) symmetric pairs of zeros $\rho =\beta +i\gamma$ and $1-\overline{\rho} = 1-\beta + i\gamma$, with $\beta \neq 1/2$. This classification will be important later. 

\section{Counting zeros in the critical strip with multiplicity}
Analytic functions have their zeros counted by the argument principle, which counts zeros with their multiplicity. Thus a zero $\rho=\beta +i\gamma$ with multiplicity $m_\rho $ is counted as $m_\rho$ zeros. To handle this situation, we replace the \lq\lq set" of zeta-zeros $Z=\{ \rho: \zeta(\rho)=0\}$ with the multiset $\mathcal{Z} = \{ \rho: \zeta(\rho)=0, \ \rho \ \text{has} \ m_\rho\ \text{copies}\}$. We now make the convention that zeros are taken from the multiset $\mathcal{Z}$. Thus, letting $ N(T)$ denote the number of zeros in the critical strip with $0< \gamma \le T$, all of the expressions here define the same number:
\begin{equation} \label{N(T)multiset} N(T) = |\{\rho\in\mathcal{Z}: 0<\gamma\le T\}| = \sum_{\substack{\rho \in \mathcal{Z}\\ 0<\gamma\le T}} 1 =\sum_{\substack{\rho \\ 0<\gamma \le T}}1 = \sum_{0<\gamma \le T} 1 .\end{equation}

We start with the Riemann-von Mangoldt formula for $N(T)$. From Titchmarsh \cite[Theorem 9.3]{Titchmarsh2} or \cite[Corollary 14.2]{MontgomeryVaughan2007} we have, for $T\ge 3$,
\begin{equation}\label{R-vM} N(T) = M(T) + \frac78+ S(T) + O(1/T),
\end{equation} 
where
\begin{equation}\label{L-S}
M(T) := \frac{T}{2\pi} \log \frac{T}{2\pi e} = \frac{1}{2\pi}\int_{2\pi e}^T\log\frac{t}{2\pi}\, dt.
\end{equation}
and von Mangoldt proved in 1905 that
\begin{equation}\label{S(T)}
S(T) := \frac1{\pi}{\rm arg }\, \zeta(\tfrac12+iT) = O(\log T).
\end{equation}
The term $O(1/T)$ is continuous, and this formula assumes $T\neq \gamma$ for any zero. If $T=\gamma$ is the ordinate of a zero, we define $N(T)= N(T^+)$ in agreement with \eqref{N(T)multiset}.\footnote[2]{The choice $N(T)= \frac12(N(T^+)+N(T^-))$ if $T$ is the ordinate of a zero is also commonly used.}
By \eqref{R-vM}, \eqref{L-S}, and \eqref{S(T)}, we have in particular,
\begin{equation} \label{N(T)sum} N(T) = \frac{T}{2\pi}\log \frac{T}{2\pi} - \frac{T}{2\pi} + O(\log{T}). \end{equation}
From \eqref{N(T)sum}, we see that
if we set
\begin{equation}\label{L} L:=\frac1{2\pi}\log T, \end{equation}
then $TL$ is the main term, and we will use this notation throughout this paper.
Thus we have
\begin{equation}\label{N(T)simple} N(T) = \frac{T}{2\pi}{\log T} + O(T) = TL + o(TL), \qquad \text{as}\qquad T\to \infty, \end{equation}
and the average vertical spacing between zeros is $1/L$.
\section{The Pair Correlation Conjecture}

Following Gallagher and Mueller \cite{GaMu78}, we now introduce the pair correlation counting function $N(T,U)$ given by
\begin{equation} \label{N(T,U)} N(T,U) := \sum_{\substack{\rho,\rho' \\ 0<\gamma, \gamma'\le T \\ 0<\gamma'-\gamma \le U}} 1.
\end{equation}
Thus we are counting the number of pairs of zeros with heights between $0$ and $T$, where the height of the second zero is greater than the first by as much as $U$. We are particularly interested when $U\asymp 1/L$. Here $f\asymp g$ means $f\ll g$ and $g\ll f$. To put this precisely, let $\lambda>0$ be a constant such that
\begin{equation} \label{lambda} UL = \lambda, \qquad \text{or}\qquad U=\frac{\lambda}{L}=\frac{2\pi\lambda}{\log T}. \end{equation} Thus $U$ has length $\lambda$ times the average spacing $1/L$. In this paper $\lambda$ is the most useful variable to consider, and therefore we define
\begin{equation} \label{N(lambda)}\mathcal{N}(\lambda) := N(T,U) =\sum_{\substack{\rho,\rho' \\ 0<\gamma, \gamma'\le T \\ 0<(\gamma'-\gamma)L \le \lambda}} 1.
\end{equation} 
Montgomery \cite{Montgomery73} made the following conjecture on pair correlation of zeros.

\begin{PCC} \label{PCC} For $UL=\lambda>0$, then
\[
\mathcal{N}(\lambda) = TL \int_0^{\lambda} \left(1 - \left(\frac{\sin\pi\alpha}{\pi \alpha}\right)^2\right) d\alpha + o(TL), \qquad \text{as}\qquad T\to \infty,
\]
uniformly in each interval $0< \lambda_0=U_0L\le \lambda=UL \le \lambda_1=U_1L<\infty$. 
\end{PCC}
\begin{remark} \label{remark1} As in \cite{GaMu78}, when applying \hyperref[PCC]{\bf PCC}, since $U_0L=\lambda_0$ and $U_1L=\lambda_1$ are arbitrary positive constants we can take $U_0$ and $U_1$ to be functions of $T$ for which $\lambda_0=U_0L\to 0$ and $\lambda_1=U_1L\to \infty$ as $T\to \infty$.
\end{remark}

\section{The Second Moment for Zeros in Short Intervals}

We now show that pair correlation of zeros is closely connected to the second moment of zeros in short intervals.
In the following proposition the second moment for zeros in short intervals is evaluated in two ways; first directly from the Riemann-von Mangoldt formula \eqref{R-vM}, and second as a double sum over pairs of zeros. 
Following Gallagher and Mueller \cite{GaMu78}, we use the forward difference notation
\begin{equation}\label{Delta} \Delta_{U} F(t):= F(t+U)-F(t). \end{equation} 

\begin{prop}[Gallagher and Mueller {\cite[Section 1]{GaMu78}}]
\label{prop1}
For $0< U \le 1$, we have as $T\to \infty$,
\begin{equation}\label{pr01eq1} \int_0^T (\Delta_UN(t))^2 \, dt = T(UL)^2 + O(TU^2L) + \int_0^T (\Delta_US(t))^2 \, dt + O(L^2), \end{equation}
and also
\begin{equation}\label{Lem3eq}
\int_0^T (\Delta_UN(t))^2 \, dt
= \sum_{\substack{\rho,\rho' \\ 0<\gamma, \gamma'\le T \\ |\gamma'-\gamma|\le U }} \left(U-|\gamma'-\gamma|\right) + O(L^2).
\end{equation}
In particular, the double sum over zeros can be rewritten as follows:
\begin{equation}\label{Lem3eq2}
\sum_{\substack{\rho,\rho' \\ 0<\gamma, \gamma'\le T \\ |\gamma'-\gamma|\le U }} \left(U-|\gamma'-\gamma|\right)
= U N^\circledast(T) + 2\int_0^UN(T,u) \,du,
\end{equation}
where 
\begin{equation}\label{N_circledast}
N^\circledast(T) := \sum_{\substack{\rho, \rho'\\ 0<\gamma,\gamma' \le T\\ \gamma = \gamma'}} 1.
\end{equation}
\end{prop}
We include the proof of this Proposition in Section \ref{sec9}. To complete the evaluation of the second moment in \eqref{pr01eq1}, we need the following unconditional result for the second moment of $\Delta_US(t)$. 

\begin{prop}[Fujii]
\label{prop2}
For $0<U\le 1$, we have as $T\to \infty$,
\begin{equation}\label{prop2eq} \int_0^T \big(\Delta_U S(t)\big)^2\, dt = \frac{T}{\pi^2} \log(2+UL) + O\left(T\sqrt{\log(2+UL)}\right). \end{equation}
\end{prop}
For the proof of this Proposition, see Fujii \cite{Fu74,Fu81} and Tsang \cite{Tsang84, Tsang86}. They both prove this unconditionally for the $2k$-th moments. This result depends on an unconditional explicit formula of Selberg for $S(t)$ \cite[Theorem 2]{Sel46} and to the authors' knowledge, there is no easy proof for that.

\begin{remark} Assuming RH, the error term in \eqref{prop2eq} can be improved to $O(T)$, but this does not improve any results we obtain from Proposition \ref{prop1}.
\end{remark}

\section{An Unconditional Average Pair Correlation Formula}

We take $UL= \lambda$ and assume $0<U\le 1$ so that $0<\lambda \le L$. Then, we immediately obtain from Propositions \ref{prop1} and \ref{prop2} that
\[ \sum_{\substack{\rho,\rho' \\ 0<\gamma, \gamma'\le T \\ |\gamma'-\gamma|L\le \lambda }} \left(\frac{\lambda}{L} -|\gamma'-\gamma|\right)
= T\left( \lambda^2 +\frac{\log(2+\lambda)}{\pi^2} + O\left(\frac{\lambda^2}{L}\right) + O\left(\sqrt{\log(2+\lambda)} \right) \right),\]
hence imposing $0<\lambda^2 \le L$ and multiplying by $L/\lambda$, we get
\begin{equation}\label{sec6_eq0}
\sum_{\substack{\rho,\rho' \\ 0<\gamma, \gamma'\le T \\ |\gamma'-\gamma|L\le \lambda }} \left(1 - \frac{|\gamma'-\gamma|L}{\lambda}\right)
= TL\left( \lambda +\frac{\log(2+\lambda)}{\pi^2\lambda} + O\left(\frac{\sqrt{\log(2+\lambda)}}{\lambda}\right) \right).
\end{equation}
By easy substitution, we see that
\[
\frac{1}{U}\int_0^UN(T,u) \,du = \frac{1}{\lambda} \int_0^\lambda \mathcal{N}(\alpha)\, d\alpha,
\]
thus dividing both sides of \eqref{Lem3eq2} by $U$, we have by \eqref{sec6_eq0},
\begin{equation} \label{averagePC} N^\circledast(T) + \frac{2}{\lambda} \int_0^\lambda \mathcal{N}(\alpha)\, d\alpha= TL\left( \lambda +\frac{\log(2+\lambda)}{\pi^2\lambda} + O\left(\frac{\sqrt{\log(2+\lambda)}}{\lambda}\right) \right).\end{equation}
This formula becomes an asymptotic formula when $\lambda \to \infty$ as $T\to \infty$.
By \eqref{averagePC} and upon using the trivial lower bound $N^\circledast(T)\ge N(T)= TL + O(T)$, we obtain unconditionally a result for repulsion of close zeros. This idea is partially implicit in Gallagher and Mueller \cite[Theorem 1]{GaMu78}.
\begin{Rep}[Gallagher and Mueller] \label{rep} For $0<\lambda\le \sqrt{L}$, we have
\begin{equation}\label{repulsion}\begin{split}
  \frac{2}{\lambda} \int_0^\lambda \Big( \mathcal{N}(\alpha) -\alpha TL \Big)\, d\alpha &= -N^\circledast(T) +TL\left( \frac{\log(2+\lambda)}{\pi^2\lambda} + O\left(\frac{\sqrt{\log(2+\lambda)}}{\lambda}\right)\right)\\ &
  \le - TL +O(T) +O\left(\frac{\log(2+\lambda)}{\lambda}TL\right) \le - TL + o_\lambda(TL), 
\end{split}\end{equation}
if $\lambda \to \infty$ as $T\to \infty$.
\end{Rep}
\begin{remark}
The Fej{\'e}r kernel factor in \hyperref[PCC]{\bf PCC} shows that close zeros tend to push each other away and thus be well-spaced and, locally, \lq\lq try to line up in approximate arithmetic progressions."
(This quote is taken from Julia Mueller's thesis \cite[p. 15]{Mueller76}.)
From \eqref{repulsion} we see this repulsion arises from the term $N^\circledast(T)$, and the negative upper bound in \eqref{repulsion} is the minimum repulsion that can occur. As we will see, \hyperref[PCC]{\bf PCC} gives this minimum repulsion. \end{remark}

\section{Horizontal multiplicity }
Recall from \eqref{N_circledast}
\[
N^\circledast(T) = \sum_{\substack{\rho, \rho'\\ 0<\gamma,\gamma' \le T\\ \gamma = \gamma'}} 1. \]
Given a zero $\rho= \beta + i\gamma$, then the solutions of the equation $\gamma' = \gamma$ are all the zeros on the horizontal line $t=\gamma$. Hence we always have the solutions $\rho =\rho'$, which we refer to as {\it diagonal terms}. In addition, if RH is false, then there are pairs of zeros with $\beta \neq 1/2$ and $\rho=\beta + i \gamma$, $\rho' = 1-\bar{\rho} = 1-\beta +i\gamma$ on the horizontal
line $t=\gamma$. We call these {\it symmetric diagonal terms}. Since we are counting pairs of zeros, if there are more than two zeros on a horizontal line, then we also get additional {\it nonsymmetric horizontal terms}. 

Soundararajan pointed out to us a more natural approach which we adopt here. Let
\begin{equation} \label{H1}
N^\circledast(T) =\sum_{\substack{\rho,\rho' \\ 0<\gamma , \gamma' \le T \\ \gamma=\gamma'}} 1 = \sum_{\substack{\rho \\ 0<\gamma \le T }}H(\gamma),
\quad \text{where} \quad
H(\gamma) := \sum_{\substack{\rho' \\ \gamma'=\gamma}} 1.
\end{equation}
If we assume RH, then $H(\gamma) = m_\rho = m_\gamma$ is the usual multiplicity of the zero $\rho$. Without RH, $H(\gamma)$ counts with multiplicity the number of zeros on the horizontal line $t=\gamma$. Thus we see, for example:

\begin{itemize}
    \item If $H(\gamma)=1$, then on the line $t=\gamma$, we have one simple zero on the critical line.
    \item If $H(\gamma)=2$, then on the line $t=\gamma$, we have either a double zero on the critical line or a pair of two symmetric simple zeros. 
    \item If $H(\gamma)=3$, then on the line $t=\gamma$, we have either a triple zero on the critical line, or a simple zero on the critical line and a pair of two symmetric simple zeros.
\end{itemize}

Note that $N^\circledast(T)$ is related to the average horizontal multiplicity for the multiset of zeros on each distinct horizontal line $t=\gamma$.
We now introduce a natural conjecture phrased using this notion of horizontal multiplicity, which is a slightly weaker version of a simplicity of zeros type conjecture.
\begin{HMH}\label{HMH} We have
\[N^\circledast(T) =\sum_{\substack{\rho \\ 0<\gamma \le T }}H(\gamma) = (1+o(1))TL \qquad \text{as} \qquad T\to \infty.\]
\end{HMH}

As it turns out, \hyperref[HMH]{\bf HMH} actually implies that asymptotically almost all the zeros are simple, and in addition, almost all of the zeros also lie on the critical line.

\begin{theorem} \label{thm2} Assuming \hyperref[HMH]{\bf HMH} is true, then asymptotically $100\%$ of the zeros of $\zeta(s)$ are simple and on the critical line.
\end{theorem}
\begin{proof}[Proof of Theorem \ref{thm2}]
If $H(\gamma)=1$ then we have exactly one zero on the horizontal line $t=\gamma$ which thus must be simple and critical. (There may be more simple and critical zeros not included here on horizontal lines with $H(\gamma)\ge 3$.) Thus, by \hyperref[HMH]{\bf HMH}, 
\[ \sum_{\substack{\rho \\ 0<\gamma \le T \\ \rho\ \text{simple and critical}}} 1 \ge \sum_{\substack{\rho \\ 0<\gamma \le T \\ H(\gamma)=1}} 1 \ge \sum_{\substack{\rho \\ 0<\gamma \le T}} (2-H(\gamma))\ge 2N(T) - (1+o(1))TL = (1+o(1))N(T). \]
\end{proof}
\begin{remark} In place of \hyperref[HMH]{\bf HMH}, we can apply this method with $N^\circledast(T) \le ({\bf C}+o(1))TL$, where $1\le \mathbf{C} <2$.  This will give us weaker proportions of the zeros depending on the size of the constant $\mathbf{C}$, see \cite[Theorem 3 and Section 7]{GS25}.
\end{remark}

\section{Proof of Theorem \ref{thm1}}
\label{secProofthm1}
By Theorem \ref{thm2}, we only need to prove that assuming \hyperref[PCC]{\bf PCC} then \hyperref[HMH]{\bf HMH} is true. 
\begin{proof}[Proof of Theorem \ref{thm1}]
Let $U_0L =\lambda_0$ and $UL =\lambda$, where $0 < \lambda_0 < 1 < \lambda$. We also take $\lambda_0\to 0$ as $T\to \infty$. If $0<u \le v$ then $0\le \mathcal{N}(u)\le \mathcal{N}(v)$, and therefore
$$ \int_0^{\lambda_0}\mathcal{N}(u)\, du \le {\lambda_0} \, \mathcal{N}(\lambda_0).
$$
By \hyperref[PCC]{\bf PCC} we have upon using the Taylor expansion of $\sin{x}$ around the origin,
$$
\mathcal{N}(\lambda_0) \ll TL\int_0^{\lambda_0} \alpha^2\, d\alpha = O(\lambda_0^3\, TL),
$$
and thus
\begin{align*} \frac{2}{\lambda}\int_0^\lambda \mathcal{N}(u) \,du &= \frac{2}{\lambda}\int_{\lambda_0}^{\lambda}\mathcal{N}(u) \,du +O\left(\frac{\lambda_0^4}{\lambda}\,TL\right)\\&
= \frac{2}{\lambda}\int_{\lambda_0}^\lambda TL\int_0^{u} \left(1 - \left(\frac{\sin\pi\alpha}{\pi \alpha}\right)^2\right)\, d\alpha \, du +o( TL)+O\left(\frac{\lambda_0^4}{\lambda}\,TL\right).
\end{align*}
We can now replace $\lambda_0$ by $0$ in the integral above with the same error $O(\lambda_0^4\,TL/\lambda)$, and on
changing the order of integration obtain
\[\begin{split}\frac{2}{\lambda}\int_0^\lambda \mathcal{N}(u) \,du &= 2TL\int_{0}^{\lambda}\left(1-\frac{\alpha}{\lambda}\right) \left(1 - \left(\frac{\sin\pi\alpha}{\pi \alpha}\right)^2\right)\, d\alpha + o( TL) +O\left(\frac{\lambda_0^4}{\lambda}\,TL\right)\\&
= \lambda TL - 2TL\int_{0}^{\lambda}\left(1-\frac{\alpha}{\lambda}\right) \left(\frac{\sin\pi\alpha}{\pi \alpha}\right)^2\, d\alpha + o(TL) +O\left(\frac{\lambda_0^4}{\lambda}\,TL\right).
\end{split}\]
By \eqref{averagePC} we conclude, since for $\lambda >1$ we have $\log (2 + \lambda) = \log \lambda + O(1/\lambda)$, 
\begin{equation} \label{endeq}\frac{N^\circledast(T)}{TL} = 2\int_{0}^{\lambda}\left(1-\frac{\alpha}{\lambda}\right) \left(\frac{\sin\pi\alpha}{\pi \alpha}\right)^2\, d\alpha +\frac{\log\lambda}{\pi^2\lambda} + o(1) + O\left( \frac{\lambda_0^4+\sqrt{\log{\lambda}}}{\lambda}\right). \end{equation}
The integral here is, for $\lambda \ge 1$,
\begin{align} \label{finalintegral}
2\int_0^{\lambda}\left(1-\frac{\alpha}{\lambda}\right) \left(\frac{\sin\pi\alpha}{\pi \alpha}\right)^2\, d\alpha = 1 - \frac{\log\lambda}{\pi^2\lambda} + O\left(\frac{1}{\lambda}\right),
\end{align}
since
\[2\int_0^{\lambda} \left(\frac{\sin\pi\alpha}{\pi \alpha}\right)^2\, d\alpha = 2\int_0^{\infty} \left(\frac{\sin\pi\alpha}{\pi \alpha}\right)^2\, d\alpha + O\left(\frac1\lambda\right) = 1 + O\left(\frac1\lambda\right), \]
and
\[\frac{2}{\lambda}\int_0^{\lambda} \alpha\left(\frac{\sin\pi\alpha}{\pi \alpha}\right)^2\, d\alpha = \frac1{\pi^2\lambda} \int_1^{\lambda} \frac{1-\cos(2\pi\alpha)}{\alpha}\, d\alpha + O\left(\frac1\lambda\right) 
= \frac{\log\lambda}{\pi^2\lambda} + O\left(\frac1\lambda\right).\]
Substituting \eqref{finalintegral} into \eqref{endeq} we conclude
\[ N^\circledast(T) = TL + o(TL)+ O\left( \frac{\lambda_0^4+\sqrt{\log{\lambda}}}{\lambda}TL\right), \qquad \text{as}\quad T\to\infty. \]

Here $\lambda_0\to 0$, and we take $\lambda\to\infty$ appropriately as $T\to \infty$ as in Remark \ref{remark1} below \hyperref[PCC]{\bf PCC}. Hence $N^\circledast(T) = TL + o(TL)$ which establishes \hyperref[HMH]{\bf HMH} and we are done.
\end{proof}

\section{Proof of the Proposition \ref{prop1}}
\label{sec9}

We will make use of the trivial estimates 
\begin{equation} \label{Delta_UN-Sbound} \Delta_UN(t)\ll (1+U)L, \quad \Delta_US(t)\ll L, \qquad \text{for} \qquad 0\le t\le T \quad \text{and} \quad 0\le U\le T, \end{equation}
which follow immediately from the estimates $S(T)\ll L$ and $N(T+1)-N(T) \ll L$ obtained from \eqref{R-vM} and \eqref{L-S}.

\begin{proof}[Proof of \eqref{Lem3eq}] This is the same proof as in \cite{GaMu78}. For $0<U\le 1$,
\[ \int_0^T (\Delta_UN(t))^2 \, dt = \int_0^T \bigg(\sum_{\substack{\rho\\ t<\gamma \le t+U}}1\bigg)^2 \, dt =
\sum_{\substack{\rho,\rho' \\ 0<\gamma, \gamma'\le T +U }}d(\gamma,\gamma'),
\]
where $d(\gamma,\gamma')$ is the measure of the set $t\in [0,T]$ with $t<\gamma,\gamma'\le t+U$,
namely,
\[
    d(\gamma, \gamma')
    :=
    \mathrm{meas}\!
    \left(
        \left[
        \max \{ \gamma - U, \gamma' - U \},
        \min \{ \gamma, \gamma' \}
        \right]
    \right)
\]
where $\mathrm{meas}$ denotes the usual Lebesgue measure on $\mathbb{R}$.
Obviously this is $0$ if $|\gamma' -\gamma|\ge U$, while if $|\gamma' -\gamma|\le U$ this measure is $= U - |\gamma'-\gamma|$ provided $0<\gamma,\gamma'\le T$. Hence, since $\gamma > 14.1$ and $0<U\le 1$, we have
\begin{align*}
\sum_{\substack{\rho,\rho' \\ 0<\gamma, \gamma'\le T +U }}d(\gamma,\gamma') 
&= 
\left(
\sum_{
\substack{
    \rho,\rho' 
    \\ 0 < \gamma , \gamma' \le T
    }
}
+
\sum_{
\substack{
    \rho,\rho' 
    \\ T < \gamma , \gamma' \le T + U
    }
}
+ \ \
2\sum_{
\substack{
    \rho,\rho' 
    \\ T < \gamma \le T +U
    \\ T-U < \gamma' \le T
    }
}
\right)
\left(U-|\gamma'-\gamma|\right) 
\\ &= \sum_{\substack{\rho,\rho' \\ 0<\gamma, \gamma'\le T\\ |\gamma'-\gamma|\le U}} \left(U- |\gamma'-\gamma|\right) + O\Bigg( \sum_{\substack{\rho,\rho' \\ T-U<\gamma, \gamma'\le T+U}}U\Bigg), 
\end{align*}
and by \eqref{Delta_UN-Sbound}, this last error term is $O(L^2)$ which proves the first equality of \eqref{Lem3eq}.
Finally, using Riemann-Stieltjes integration and recalling the definition of $N^\circledast(T)$ in \eqref{N_circledast}, we have
\[\begin{split} \sum_{\substack{\rho,\rho' \\ 0<\gamma, \gamma'\le T\\ |\gamma'-\gamma|\le U}} \left(U-|\gamma'-\gamma|\right) &= U N^\circledast(T) + 2\sum_{\substack{\rho,\rho' \\ 0<\gamma, \gamma'\le T \\ 0< \gamma'-\gamma\le U } }\left(U-(\gamma'-\gamma)\right) \\& =U N^\circledast(T) + 2\int_0^U(U- u)\, d_uN(T,u) \\&= U N^\circledast(T) + 2\int_0^UN(T,u) \,du. \end{split}\]
\end{proof}

\begin{proof}[Proof of \eqref{pr01eq1}]

Let $0< U\le 1$. 
By \eqref{R-vM}, we have
\begin{equation} \label{DeltaN}
\begin{split}
\int_0^T (\Delta_UN(t))^2 \, dt = \int_2^T (\Delta_UN(t))^2 \, dt
&= \int_2^T\left(\Delta_UM(t) + \Delta_US(t) + O\left(\tfrac1t\right)\right)^2 \, dt
\\
&= \int_2^T\left(\Delta_UM(t) \right)^2\, dt + \int_2^T \left(\Delta_US(t) + O\left(\tfrac1t\right)\right)^2\, dt \\
&\qquad\qquad+ O\left(\left|\int_2^T \Delta_UM(t)\left(\Delta_US(t) + O\left(\tfrac1t\right) \right)\,dt\right|\right).
\end{split}\end{equation}
Now
\[\Delta_UM(t)=\frac{1}{2\pi}\int_{t}^{t+U}\log\frac{u}{2\pi}\, du = \frac{U}{2\pi}(\log t + O(1)),\]
and therefore
\begin{equation}\label{DeltaL}
\int_2^T\left(\Delta_UM(t) \right)^2\, dt = T(UL)^2 + O(T U^2 L). \end{equation}
Since $\Delta_US(t) \ll L$, 
\begin{equation}\begin{split} \label{DeltaS}
\int_2^T\left( \Delta_US(t) + O\left(\tfrac1t\right) \right)^2\, dt &= \int_2^T\left( \Delta_US(t)\right)^2\, dt + O\left( \int_2^T \frac{|\Delta_US(t)|}t\, dt\right) + O(1) \\
&= \int_2^T \left(\Delta_US(t)\right)^2\, dt + O\left(L\int_2^T \frac1t\, dt\right) + O(1) \\
&= \int_0^T \left(\Delta_US(t)\right)^2\, dt + O\left(L^2\right). \end{split}\end{equation}

Finally, 
\begin{equation}\begin{split}\label{DeltaLS}
\int_2^T \Delta_UM(t)\left(\Delta_US(t) + O\left(\tfrac1t\right)\right)\, dt
&= \int_2^T \Delta_UM(t)\Delta_US(t)\,dt + O\left(\int_2^T \Delta_UM(t)\,\frac{dt}{t} \right)\\
&= \int_2^T \Delta_UM(t)\Delta_US(t)\,dt + O(UL^2). \end{split}\end{equation}
For $2\le t\le T$, we have
\[\begin{split}
S_1(t,U) :&= \int_2^t \Delta_US(u)\,du = \int_2^t S(u+U)\,du - \int_2^tS(u)\,du \\
&= \int_t^{t+U}S(u)\,du 
- \int_2^{2+U} S(u)\,du \ll U\log t. \end{split}\]
Since $\Delta_UM(t)$ is monotonically increasing 
with $\frac{d}{dt}\Delta_UM(t) = \frac{1}{2\pi} \log(1 + \frac{U}{t})$, we have 
\begin{equation}\begin{split}\label{DeltaLS2}\int_2^T \Delta_UM(t)\Delta_US(t)\, dt &= \int_2^T \Delta_UM(t)\frac{d}{dt} (S_1(t,U))\, dt\\& 
= \Delta_UM(t) S_1(t,U) \bigg|_2^T -\frac{1}{2\pi}\int_2^T \log\!\left(1 + \frac{U}{t}\right)S_1(t,U)\, dt \\&
\ll U^2 L^2 + U^2 \int_2^T\frac{\log t}{t} \, dt
\ll (UL)^2 \ll L^2. \end{split}\end{equation}
Combining \eqref{DeltaN}, \eqref{DeltaL}, \eqref{DeltaS}, \eqref{DeltaLS}, and \eqref{DeltaLS2}, we have
\[ \int_0^T (\Delta_UN(t))^2 \, dt = T(UL)^2 + O(TU^2L)+ \int_1^T (\Delta_US(t))^2 \, dt + O(L^2). \]
\end{proof}

\section*{Acknowledgement of Funding}

DAG, JL and AIS thank the American Institute of Mathematics for its hospitality and for providing a pleasant research environment where this research project first started. JL was supported by the National Research Foundation of Korea (NRF) grant funded by the government of the Republic of Korea (MSIT) (No. RS-2025-00553763, RS-2024-00415601, and RS-2024-00341372). JS was supported by an AMS-Simons Research Enhancement Grant for Primarily Undergraduate Institution (PUI) Faculty. AIS was supported by the Inamori Research Grant 2024, JSPS KAKENHI Grant Number 22K13895, and Kyushu University International Research Leader Training Program (EBXU0101).

\section*{Additional Statements}
On behalf of all authors, the corresponding author states that there is no conflict of interest.
Data sharing is not applicable to this article as no datasets were generated or analyzed during the current study.

\bibliographystyle{alpha}
\bibliography{AHReferences}

\end{document}